\documentclass[a4paper]{amsart}

\usepackage{amsthm,amsfonts,amsmath, amssymb}

\usepackage{isomath}
\usepackage{lmodern} 
\usepackage[italic]{mathastext}

\newcommand{\cvec}[1]{\begin{pmatrix}#1\end{pmatrix}} 
\let\oldmarginpar\marginpar
\renewcommand{\marginpar}[1]{\oldmarginpar{\small\textit{{#1}}}}
\newcommand{\mpar}[1]{\marginpar{#1}}

\usepackage{setspace}
\usepackage{graphicx,float,enumerate,subfigure,tikz}
\usepackage[ruled,boxed,commentsnumbered,norelsize]{algorithm2e}
\usepackage{verbatim}
\usepackage{url}   
\usepackage[capitalise]{cleveref}

\usepackage[sort&compress,comma]{natbib}
 \usepackage{enumitem}
\setlist[enumerate,1]{label=(\roman*),font=\normalfont}

\newtheorem{theorem}{Theorem}
\newtheorem{corollary}[theorem]{Corollary}
\newtheorem{lemma}[theorem]{Lemma}
\newtheorem{proposition}[theorem]{Proposition}

\newenvironment{claimproof}{\noindent\textit{Proof.}}{\hfill$\square$}

\newtheorem{conjecture}[theorem]{Conjecture}

\theoremstyle{definition}

\theoremstyle{remark}
\newtheorem{remark}[theorem]{Remark}

\crefname{remark}{Remark}{Remarks}

\theoremstyle{definition}
\newtheorem{claim}{Claim}
\usepackage{etoolbox}
\AtEndEnvironment{proof}{\setcounter{claim}{0}}

\theoremstyle{remark}
\newtheorem{case}{Case}
\usepackage{etoolbox}
\AtEndEnvironment{proof}{\setcounter{case}{0}}

\crefname{rmk}{Remark}{Remarks} 
\crefname{problem}{Problem}{Problems}

\usepackage{mathtools}
\DeclarePairedDelimiter\ceil{\lceil}{\rceil}
\DeclarePairedDelimiter\floor{\lfloor}{\rfloor}

\newcommand{\R}{\mathbb{R}}

\renewcommand{\ge}{\geqslant} 
\renewcommand{\le}{\leqslant}

\date{\today}
\title{Linkedness of Cartesian products of complete graphs}

\usepackage{amsaddr}

\author{Leif K.~ J{\o}rgensen}
\address{Department of Mathematical Sciences, Aalborg University, Denmark}
\email{\texttt{leifkjorgensen@gmail.com}}
\author{Guillermo Pineda-Villavicencio}
\address{Centre for Informatics and Applied Optimisation, Federation University,  Australia\\School of Information Technology, Deakin University, Geelong, Australia}  
\email{\texttt{work@guillermo.com.au}}
\author{Julien Ugon}
\address{Centre for Informatics and Applied Optimisation, Federation University,  Australia\\School of Information Technology, Deakin University, Geelong, Australia}
\email{\texttt{julien.ugon@deakin.edu.au}}

\thanks{Thanks}
 
\begin{document}
\begin{abstract} This paper is concerned with the linkedness of Cartesian products of complete graphs. A graph with at least $2k$ vertices is {\it $k$-linked} if, for every set of $2k$ distinct vertices organised in arbitrary $k$ pairs of vertices, there are $k$ vertex-disjoint paths joining the vertices in the pairs. 

We show that the Cartesian product $K^{d_{1}+1}\times K^{d_{2}+1}$ of complete graphs $K^{d_{1}+1}$ and $K^{d_{2}+1}$ is $\floor{(d_{1}+d_{2})/2}$-linked for $d_{1},d_{2}\ge 2$, and this is best possible.
 

This result is connected to graphs of simple polytopes. The Cartesian product $K^{d_{1}+1}\times K^{d_{2}+1}$  is the graph of the Cartesian product $T(d_{1})\times T(d_{2})$ of a $d_{1}$-dimensional simplex $T(d_{1})$ and a $d_{2}$-dimensional simplex $T(d_{2})$. And the polytope $T(d_{1})\times T(d_{2})$ is a {\it simple polytope}, a $(d_{1}+d_{2})$-dimensional polytope in which every vertex is incident to exactly $d_{1}+d_{2}$ edges. 
 
While not every $d$-polytope is $\floor{d/2}$-linked, it may be conjectured that every simple $d$-polytope is.  Our result implies the veracity of the revised conjecture for Cartesian products of two simplices. \end{abstract}
\maketitle  

\section{Introduction}

Denote by $V(X)$ the vertex set of a graph. Given sets $A,B$ of vertices in a graph, a path from $A$ to $B$, called an {\it $A-B$ path}\mpar{$A-B$ path}, is a (vertex-edge) path $L:=u_{0}\ldots u_{n}$ in the graph such that $V(L)\cap A=\{u_{0}\}$  and $V(L)\cap B=\{u_{n}\}$. We write $a-B$ path instead of $\{a\}-B$ path, and likewise, write $A-b$ path instead of $A-\{b\}$.

Let $G$ be a graph and $X$ a subset of $2k$ distinct vertices of $G$. The elements of $X$ are called {\it terminals}\mpar{terminals}. Let $Y:=\{\{s_{1},t_{1}\}, \ldots,\{s_{k},t_{k}\}\}$ be an arbitrary labelling and (unordered) pairing of all the vertices in $X$. We say that $Y$ is {\it linked}\mpar{linked pairs} in $G$ if we can find disjoint $s_{i}-t_{i}$ paths for $i\in [1,k]$, the interval $1,\ldots,k$. The set $X$ is {\it linked}\mpar{linked set} in $G$ if every such pairing of its vertices is linked in $G$. Throughout this paper, by a set of disjoint paths, we mean a set of vertex-disjoint paths. If $G$ has at least $2k$ vertices and every set of exactly $2k$ vertices is linked in $G$, we say that $G$ is {\it $k$-linked}\mpar{$k$-linked}. 

This paper studies the linkedness of Cartesian products of complete graphs. Linkedness of Cartesian products has been studied in the past \citep{Mes16}. The {\it Cartesian product} $G_{1}\times G_{2}$ of two graphs $G_{1}$ and $G_{2}$ is the graph defined on the pairs $(v_{1},v_{2})$ with $v_{i}\in G_{i}$ and with two pairs $(u_{1},u_{2})$ and $(v_{1},v_{2})$ being adjacent if, for some $\ell\in \{1,2\}$, $u_{\ell}v_{\ell}\in E(G_{\ell})$  and $u_{i}=v_{i}$ for $i\ne \ell$. We prove that the Cartesian product $K^{d_{1}+1}\times K^{d_{2}+1}$ of complete graphs $K^{d_{1}+1}$ and $K^{d_{2}+1}$ is $\floor{(d_{1}+d_{2})/2}$-linked for $d_{1},d_{2}\ge 0$, and that there are products that are not $\floor{(d_{1}+d_{2}+1)/2}$-linked; hence this result is best possible. Here $K^{t}$ denotes the complete graph on $t$ vertices.

Our result is connected to questions on the linkedness of a polytope. A (convex) polytope is the convex hull of a finite set $X$ of points in $\R^{d}$; the \textit{convex hull} of $X$ is  the smallest convex set containing $X$.  The \textit{dimension} of a polytope in $\R^{d}$ is one less than the maximum number of affinely independent points in the polytope; a set of points $\vec p_{1},\ldots, \vec p_{k}$ in $\R^{d}$ is {\it affinely independent} if  the $k-1$ vectors $\vec p_{1}-\vec p_{k},\ldots, \vec p_{k-1}-\vec p_{k}$ are linearly independent.  A polytope of dimension $d$ is referred to as a \textit{$d$-polytope}.

 The {\it Cartesian product} $P\times P'$ of a $d$-polytope $P\subset \R^{d}$ and a $d'$-polytope $P'\subset \R^{d'}$ is the Cartesian product of the sets $P$ and $P'$: 
\begin{equation*}
P\times P'=\left\{\cvec{p\\p'}\in \R^{d+d'}\middle|\; p\in P,\, p'\in P\right\}.
\end{equation*} 
The resulting polytope is $(d+d')$-dimensional. The {\it graph} $G(P)$ of a polytope $P$ is the undirected graph formed by the vertices and edges of the polytope. It follows that the graph $G(P\times P')$ of the Cartesian product $P\times P'$ is the Cartesian product $G(P)\times G(P')$ of the graphs $G(P)$ and $G(P')$. 

A $d$-simplex $T(d)$ is the convex hull of $d+1$ affinely independent points in $\R^{d}$. The graph of $T(d)$ is the complete graph $K^{d+1}$. As a consequence, our result implies that the graph of the Cartesian product $T({d_{1}})\times T({d_{2}})$ is $\floor{(d_{1}+d_{2})/2}$-linked for $d_{1},d_{2}\ge 0$. Henceforth, if the graph of a polytope is $k$-linked we say that the polytope is also {\it $k$-linked}.

The first edition of the Handbook of Discrete and Computational Geometry \cite[Problem 17.2.6]{GooORo97-1st} posed the question of whether or not every $d$-polytope is $\floor{d/2}$-linked. This question was answered in the negative by \cite{Gal85}. None of the known counterexamples are {\it simple $d$-polytopes}, $d$-polytopes in which every vertex is incident to exactly $d$ edges. Hence, it may be hypothesised that the conjecture holds for such polytopes. 

\begin{conjecture}Every simple $d$-polytope is $\floor{d/2}$-linked for $d\ge 2$.  
\end{conjecture}

Cartesian products of simplices are simple polytopes, and so our result supports this revised conjecture. Furthermore, Cartesian products of simplices and duals of cyclic polytopes are related; the dual of a cyclic $d$-polytope with $d+2$ vertices is the Cartesian product of a ${\floor{d/2}}$-simplex and  a ${\ceil{d/2}}$-simplex \citep[Ex.~0.6]{Zie95}. Hence we obtain that the dual of a cyclic $d$-polytope on $d+2$ vertices is also  $\floor{d/2}$-linked for $d\ge 2$.   

Unless otherwise stated, the graph theoretical notation and terminology follows from \cite{Die05} and the polytope theoretical notation and terminology from \cite{Zie95}. Moreover, when referring to graph-theoretical properties of a polytope such as linkedness and connectivity, we mean properties of its graph.

\section{Linkedness of Cartesian products of complex graphs}
  \label{sec:simplex-products}
  
The contribution of this section is a sharp theorem (\cref{thm:simplex-products}) that tells the story of the linkedness of Cartesian product of two complete graphs.

\begin{figure}  
\includegraphics{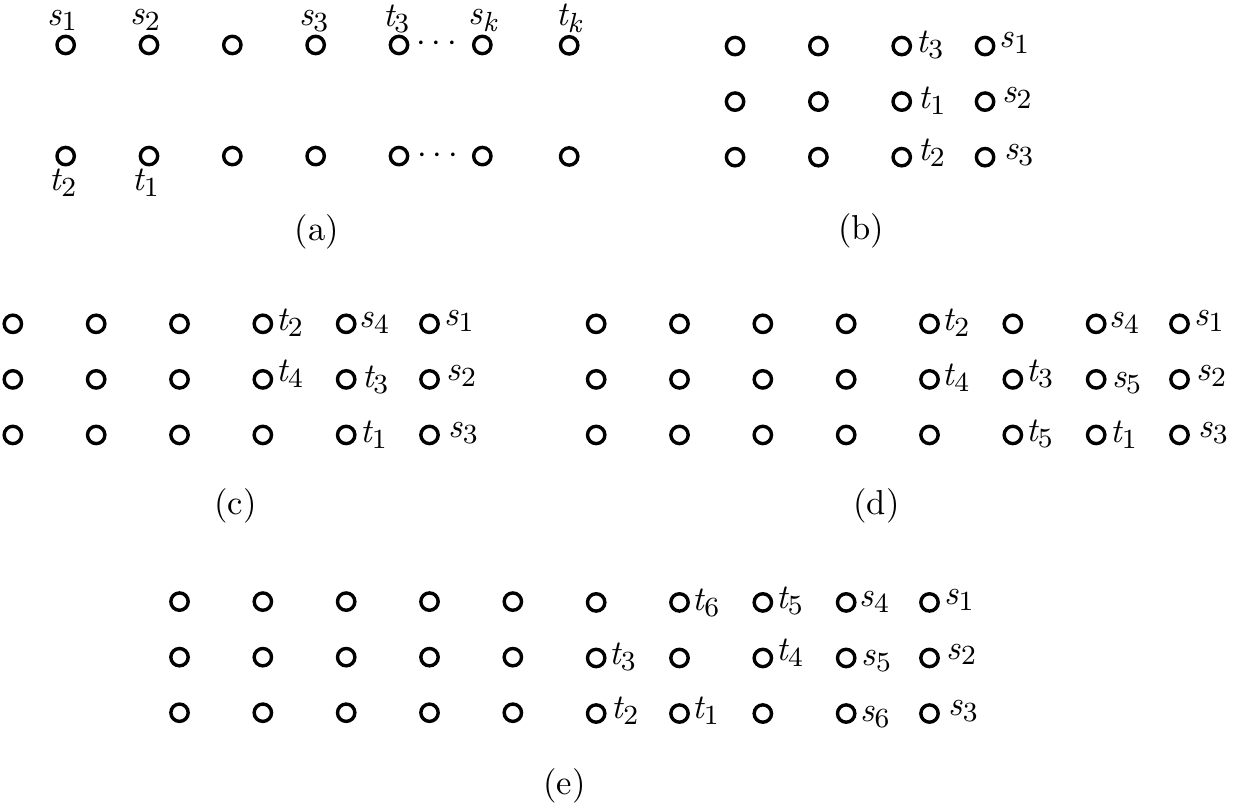}
\caption{No feasible linkage problems for $K^{d_{1}+1}\times K^{d_{2}+1}$, $k=\floor{(d_{1}+d_{2}+1)/2}$, $d_{1}\le 2$ and $d_{2}> d_{1}$. (a) The case $d_{1}=1$ and $d_{2}> d_{1}$. (b) The case $d_{1}=2$ and $d_{2}=3$. (c) The case $d_{1}=2$ and $d_{2}=5$. (d) The case $d_{1}=2$ and $d_{2}=7$. (e) The case $d_{1}=2$ and $d_{2}=9$. Each row of each part (a)-(e) is a complete graph whose edges have not been drawn.}\label{fig:no-link}
\end{figure}

\begin{theorem}
\label{thm:simplex-products} The Cartesian product of two complete graphs  $K^{d_{1}+1}$ and $K^{d_{2}+1}$ is $\floor{(d_{1}+d_{2})/2}$-linked for every $d_{1},d_{2}\ge 0$.\end{theorem}

\begin{remark}\label{rmk:sharp-simplex-products}\cref{thm:simplex-products} is best possible. There are products $K^{d_{1}+1}\times K^{d_{2}+1}$ 
that are not $\floor{(d_{1}+d_{2}+1)/2}$-linked:
\begin{enumerate}
\item $K^{2}\times K^{d_{2}+1}$ for even $d_{2}\ge 1$, and
\item $K^{3}\times K^{d_{2}+1}$ for $d_{2}=1,3,5,7,9$. 
\end{enumerate}   
For each of these cases, \cref{fig:no-link} provides a pairing of terminals  that cannot be $\floor{(d_{1}+d_{2}+1)/2}$-linked. We conjecture these are the only such cases. 
\end{remark}

An immediate corollary of \cref{thm:simplex-products} is the following.  
\begin{corollary}
\label{cor:simplex-products} The Cartesian product of two simplices $T({d_{1}})$ and $T({d_{2}})$ is $\floor{(d_{1}+d_{2})/2}$-linked for every $d_{1},d_{2}\ge 0$.\end{corollary}

The notions of linkage, linkage problem, and valid path will simplify our arguments. A {\it linkage} in a graph is a subgraph in which every component is a path. Let  $X$ be a set of vertices in a graph and let $Y:=\{\left\{s_{1},t_{1}\right\},\ldots,\{s_{k},t_{k}\}\}$ be a  pairing of all the vertices of $X$. A {\it $Y$-linkage} $\{L_{1},\ldots,L_{k}\}$ is a set of disjoint paths with the path $L_{i}$ joining the pair $\{s_{i},t_{i}\}$ for $i=1,\ldots,k$.  We may also say that $Y$ represents our {\it linkage problem}, and if $Y$ is linked in $G$ then our linkage problem is {\it feasible} and {\it infeasible} otherwise. A path in the graph is called {\it $X$-valid} if no inner vertex of the path is in $X$. 
Let  $X$ be a set of vertices in a graph $G$. Denote by $G[X]$ the subgraph of $G$ induced by $X$, the subgraph of $G$ that contains all the edges of $G$ with vertices in $X$. Write $G-X$ for $G[V(G)\setminus X]$. 
 
Consider a linkage problem $Y:=\{\{s_{1},t_{1}\}, \ldots,\{s_{k},t_{k}\}\}$ on a set $X$ of $2k$ vertices in a graph $G$. Consider a linkage $\mathcal L$ from a subset $Z$ of $X$ to some set $Z'$ disjoint from $X$ and label the vertices of $Z'$ such that the path in $\mathcal L$ with end $z_{i}\in Z$ has its other end $z_{i}'\in Z'$. Then the linkage $\mathcal L$ in $G$ {\it induces} a  linkage problem $Y'$ in $(G-V(\mathcal L))\cup Z'$ where the vertices of $X\setminus Z$  remain and the vertices of $Z$ have been replaced by the vertices of $Z'$. Slightly abusing terminology, we also call terminals the vertices of $Z'$. If the problem $Y'$ is feasible in $(G-V(\mathcal L))\cup Z'$, so is the problem $Y$ in $G$.

Since we make heavy use of Menger's theorem \cite[Thm.~3.3.1]{Die05}, we next remind the reader of one of its consequences.

\begin{theorem}[Menger's theorem]\label{thm:Menger} Let $G$ be a $k$-connected graph, and let $A$ and $B$ be two subsets of its vertices, each of cardinality at least $k$. Then there are $k$ disjoint $A-B$ paths in $G$. 
\end{theorem}

We fix some notation and terminology for the remaining of the section. Let $G$ denote the graph $K^{d_{1}+1}\times K^{d_{2}+1}$. We think of $G=K^{d_{1}+1}\times K^{d_{2}+1}$ as a grid with $d_{1}+1$ rows and $d_{2}+1$ columns. In this way, the entry in Row $i$ and Column $j$ can be referred to as $G[i,j]$.
 
When we write about a row $r$ of  subgraph $G'$ of $G$, we think of $r$ as a subgraph of $G'$ and as the number $r$ so that we can write about the $r$th row of $G'$ or $G$; this ambiguity should cause no confusion. An entry in the grid $K^{d_{1}+1}\times K^{d_{2}+1}$ with no terminal is said to be {\it free}, as is a row or a column of a subgraph of $G$ with no terminal. A row or a column of a subgraph of $G$ with every entry being occupied by a terminal is said to be {\it full}.  

\begin{figure}   
\includegraphics{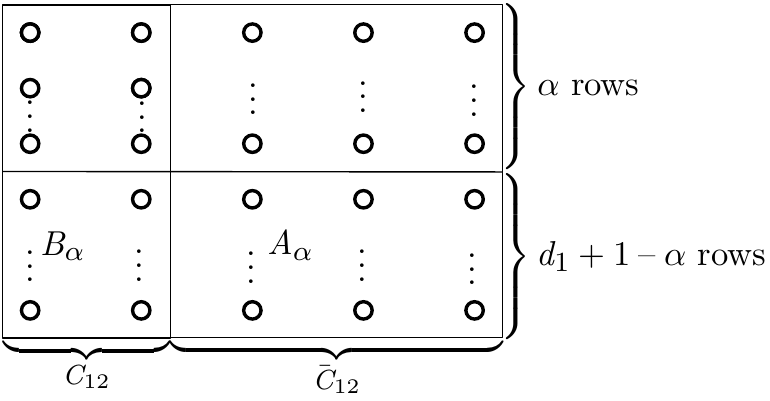}
\caption{Depiction of the subgraphs $B_{\alpha+1}$, $A_{\alpha+1}$, $C_{12}$, and $\bar C_{12}$ of $K^{d_{1}+1}\times K^{d_{2}+1}$.}\label{fig:aux-fig-products}
\end{figure}

We need the following induced subgraphs of  $G$: 
\begin{align*}
C_{ab\ldots z}&, \text{the subgraph formed by the union of Columns $a,b,\ldots,z$;} \\  
\bar C_{ab\ldots z}&, \text{the subgraph obtained by removing Columns $a,b,\ldots,z$;}\\
R_{ab\ldots z}&, \text{the subgraph formed by the union of Rows $a,b,\ldots,z$;} \\  
\bar R_{ab\ldots z}&, \text{the subgraph obtained by removing Rows $a,b,\ldots,z$;}\\
A_{\alpha}&, \text{the induced subgraph of  $\bar C_{12}$ obtained by removing its first $\alpha$ rows; and}\\
B_{\alpha}&, \text{the subgraph of $C_{12}$ obtained by removing its first $\alpha$ rows.}
\end{align*}
For instance, $\bar C_{1}$ denotes the subgraph of $G$ obtained by removing the first column,  $C_{12}$  the subgraph formed by the first two columns of $G$,  and $\bar C_{12}$ denotes the subgraph obtained by removing the first two columns of $G$; observe $\bar C_{12}$ is isomorphic to $K^{d_{1}+1}\times K^{d_{2}-1}$.  \Cref{fig:aux-fig-products} depicts some of the aforementioned subgraphs of $K^{d_{1}+1}\times K^{d_{2}+1}$. 

The connectivity of $K^{d_{1}+1}\times K^{d_{2}+1}$ is stated below. 

\begin{lemma}[{\citealt[Thm.~1]{Spa08}}]   
\label{lem:connectivity-complete-graphs} The (vertex)connectivity of $K^{d_{1}+1}\times K^{d_{2}+1}$ is precisely $d_{1}+d_{2}$.
\end{lemma}

We continue fixing further notation. Henceforth let $k:=\floor{(d_{1}+d_{2})/2}$. And let $X$ be a subset of $2k$ vertices of $G$ and let $Y:=\{\{s_{1},t_{1}\}, \ldots,\{s_{k},t_{k}\}\}$ be a pairing of all the vertices in $X$.

We first settle the simple cases of $(0,d_{2})$ and $(1,d_2)$ for $d_{2}\ge 0$. 

\begin{proposition}[Base cases]\label{prop:basic} For $d_{2}\ge 0 $ the Cartesian products $K^{1}\times K^{d_{2}+1}$ and $K^{2}\times K^{d_{2}+1}$ are both  $\floor{(1+d_{2})/2}$-linked. This statement is best possible.
\end{proposition}

\begin{proof}The lemma is true for the pair $(0,d_{2})$  for each $d_{2}\ge 0$, since $K^{1}\times K^{d_{2}+1}=K^{d_{2}+1}$ and $K^{d_{2}+1}$ is $\floor{(1+d_{2})/2}$-linked. This is best possible. 

 The graph $K^{2}\times K^{d_{2}+1}$ is $(1+d_{2})$-connected by \cref{lem:connectivity-complete-graphs}.  Use Menger's theorem~\eqref{thm:Menger} to bring the $1+d_{2}$ terminals  to the subgraph $\bar R_{1}$ through a linkage $\{S_{1},\ldots,S_{k},T_{1},\ldots,T_{k}\}$ with $S_{i}:=s_{i}-\bar R_{1}$ and $T_{i}:=t_{i}-\bar R_{1}$ for $i\in[1,k]$. Letting $\{\bar s_{i}\}:=V(S_{i})\cap V(\bar R_{1})$ and $\{\bar t_{i}\}:=V(T_{i})\cap V(\bar R_{1})$, we produce a new linkage problem $Y':=\{(\bar s_{1},\bar t_{1}\},\ldots,(\bar s_{k},\bar t_{k}\}\}$ in $\bar R_{1}$ whose feasibility implies that of $Y$ in $G$.  To solve $Y'$ link the pairs of $Y'$  in the subgraph $\bar R_{1}$, which is isomorphic to $K^{d_{2}+1}$, using the $\floor{(1+d_{2})/2}$-linkedness of $K^{d_{2}+1}$. \Cref{fig:no-link}(a) shows an infeasible linkage problem with $\floor{(2+d_{2})/2}$ pairs in the graph $K^{2}\times K^{d_{2}+1}$.
\end{proof}


In what follows we aim to find a $Y$-linkage $\{L_{1},\ldots,L_{k}\}$ in $G$ with $L_{i}$ joining the pair $\{s_{i},t_{i}\}$ of $Y$ for $i\in [1,k]$.  Our proof is by induction on $(d_{1},d_{2})$ with the base cases settled in \cref{prop:basic}. If there is a pair of $Y$, say $\{s_{1},t_{1}\}$, lying in some column or row of $G$, say in Column 1, we send every terminal $s_{i}\in C_{1}$ that is different from $s_{1}$ and $t_{1}$ and that is not adjacent to $t_{i}$ to the subgraph $\bar C_{1}$, and apply the induction hypothesis on $\bar C_{1}$. Otherwise, we may assume every pair of $Y$ lies in two distinct columns or rows, say the pair $\{s_{1},t_{1}\}$ lies in $C_{12}$; then we send every terminal $s_{i}\in C_{12}$ that is different from $s_{1}$ and $t_{1}$ and that is not adjacent to $t_{i}$ to the subgraph $\bar C_{12}$, and apply the induction hypothesis to $\bar C_{12}$. We develop these ideas below.
 
The definition of $k$-linkedness gives the following lemma at once; we will use it implicitly hereafter.

\begin{lemma}\label{lem:k-linked-def} Let $\ell\le k$. Let $X$ be a set of $2\ell$ distinct vertices of a $k$-linked graph $K$,  let $Y$ be a labelling and pairing of the vertices in $X$, and let $Z$ be a set of $2k-2\ell$ vertices in $K$ such that $X\cap Z=\emptyset$. Then there exists a $Y$-linkage in $K$ that avoids every vertex in $Z$.  
\end{lemma}

Besides, basic algebraic manipulation yields the following inequality.

\begin{lemma}\label{lem:inequality}  
If $x\ge2$ and $y\ge 2$ then $x(y-1)>x+y-3$.
\end{lemma}
\begin{proof}
The inequality simplifies to $(x-1)(y-2)>-1$.
\end{proof}
 
We are now ready to put together all the elements of the proof of \cref{thm:simplex-products}. 

\begin{proof}[Proof of \cref{thm:simplex-products}] Let $k:=\floor{(d_{1}+d_{2})/2}$. Then $d_{1}+d_{2}\ge 2k$. 

\cref{prop:basic} gives the result for the pairs $(d_{1}, 0)$, $(0,d_{2})$, $(d_{1}, 1)$, and $(1,d_{2})$  for each $d_{1},d_{2}\ge 0$. Hence,  our bidimensional induction on $(d_{1}, d_{2})$ can start with the assumption of $d_{1},d_{2}\ge 2$.  

We first deal with the case where a pair in $Y$, say $\{s_{1},t_{1}\}$, lies in some column or some row of $G$, say in Column 1. 

\begin{case} A pair in $Y$, say $\{s_{1},t_{1}\}$, lies in Column 1.
\end{case}            

The induction hypothesis ensures that the subgraph $\bar C_{1}$ is $(k-1)$-linked. Hence it suffices to show that all the terminals in $C_1$ other than $s_1,t_1$ can be moved to $\bar C_1$ via a linkage; Menger's theorem \eqref{thm:Menger} guarantees this. 

 Let $U$ be the set of terminals in $C_{1}$ other than $s_{1}$ and $t_{1}$, and let $W$ be the set of terminals in $\bar C_{1}$. Then $|U|+|W|=d_1+d_2-2$. Besides, the subgraph $G-(W\cup \{s_1,t_1\})$ is $|U|$-connected, as $G$ is $(d_1+d_2)$-connected (\cref{lem:connectivity-complete-graphs}).  In the case of $d_1,d_2\ge 2$, \cref{lem:inequality} yields that   $\bar C_1$ has more than $|U\cup W|$ vertices:
\begin{align*}
   |\bar C_1|=(d_1+1)d_2>d_1+1+d_2+1-3>d_1+d_2-2=|U|+|W|.
\end{align*}
Menger's theorem \eqref{thm:Menger} applies and gives disjoint $C_{1}-\bar C_{1}$ paths from the terminals in $U$ to $|U|$ free entries in $\bar C_1$. The $(k-1)$-linkedness of $\bar C_{1}$ now settles the case.
   	
By symmetry, we can assume that every pair $\{s_{i},t_{i}\}$ in $Y$ lies in two different columns or rows and that $s_{i},t_{i}$ are not adjacent.  Without loss of generality, assume that
 \begin{equation}\label{eq:pair-C12}
\parbox{0.8\textwidth}{\it  $s_{1}$ is Column 1  and $t_{1}$ is in Column 2 of $C_{12}$.} \tag{*} 
\end{equation}  
 
The induction hypothesis also ensures that  both $\bar C_{12}$ and $\bar R_{12}$ are $(k-1)$-linked.  We consider two further cases based on the number of terminals  in $C_{12}$ or $R_{12}$. 

\begin{case}\label{case:d1+2}  The subgraph $C_{12}$ contains precisely $d_{1}+2-\alpha$ terminals, including $\{s_{1},t_{1}\}$, where $0\le \alpha\le d_{1}$. 
\end{case} 

Excluding $\{s_{1},t_{1}\}$, there are at most $d_{1}$ terminals in $C_{12}$, and there are $d_{1}+1$  internally-disjoint $s_{1}-t_{1}$ paths in $C_{12}$ of length at most three: two length-two paths and $d_{1}-1$ length-three paths. One of these $s_{1}-t_{1}$ paths, say $L_{1}$,  avoids every other terminal in $C_{12}$. 

Without loss of generality, assume that Row 1 in $C_{12}$ is part of the path  $L_{1}$; that is, 
\begin{equation}\label{eq:path-L1}
\left\{G[1,1],G[1,2]\right\} \subseteq V(L_{1}).\tag{**} 
\end{equation}  
It is true that $(V(L_{1})\cap V(B_{1}))\subseteq \left\{s_{1}, t_{1}\right\}$.    

In the subcase $\alpha=d_{1}$, every pair in $Y\setminus \{s_1,t_1\}$ is in $\bar C_{12}$, and the induction hypothesis on $\bar C_{12}$ settles the subcase. 

Suppose that $\alpha=d_{1}-1$, say  $C_{12}$ contains $\{s_{1},t_{1},s_{2}\}$. Then $s_2\in B_1$ and $t_2\in \bar C_{12}$. We may assume $s_1,s_2$ are in Column 1 and $t_1$ is in Column  2. We show there is an $X$-valid $s_2-A_1$ path $L_2'$ such that the vertex $x\in V(L_2')\cap V(A_1)$ is either $t_2$ or a nonterminal.

Through each entry of Column 1 of $B_1$, there is a $s_2-A_1$ path of length at most two: one of length one and $d_1-1$ of length two. Moreover, $d_1-1$ of such paths avoid $s_1$. To ensure the existence of $L_2'$, it suffices to show that $A_1$ cannot have $d_1-1$ rows that are full of terminals other than $t_2$. According to \cref{lem:inequality}, the inequality 
\begin{align*}
(d_{1}-1)(d_{2}-1)>d_{1}-1+d_{2}-3=| X\setminus\{s_1,t_1,s_2,t_2\}|
\end{align*}
holds for $d_{1},d_{2}\ge 2$. Hence we get the existence of $L_2'$. As a result,  the solution of the new problem $Y':=\left\{\left\{x,t_{2}\right\},\left\{s_{3},t_{3}\right\},\ldots, \left\{s_{k},t_{k}\right\}\right\}$ in $\bar C_{12}$ induces a solution of the problem $Y$ in $G$.  And the solution of $Y'$ follows from the $(k-1)$-linkedness of $\bar C_{12}$. 

Henceforth assume that $\alpha\le d_1-2$. To finalise \cref{case:d1+2},  we require a couple of claims. 
 
\begin{claim}
Suppose that there are at most $d_1+2$ terminals in $B_1=K^{d_1}\times K^2$. Then there is an injection from the set of rows of $B_{1}$ that contain two terminals $x_{1},x_{2}$ such that $\left\{x_{1},x_{2}\right\}\cap \left\{s_1,t_{1}\right\}=\emptyset$ to the set of rows of $B_{1}$ that contain either no terminal or a terminal in $\{s_1,t_{1}\}$ but no other terminal. 
\label{cl:full-free-rows}
\end{claim}  
	
\begin{claimproof}
This follows from a simple counting argument. The number of rows in $B_{1}$ is $d_{1}$. Let $m$ denote the number of rows of $B_{1}$ that contain two terminals $x_{1},x_{2}$ such that $\left\{x_{1},x_{2}\right\}\cap \left\{s_1,t_{1}\right\}=\emptyset$ and let $n:=|(X\cap V(B_{1}))\setminus \left\{s_1,t_{1}\right\}|$. It follows that the number of rows of $B_{1}$ that contains precisely one terminal $x\not\in \{s_1,t_{1}\}$ is $n-2m$; either $s_1$ or $t_{1}$ may be in these rows. As a result, the number of rows of $B_{1}$ that contain either no terminal or a terminal in $\{s_1,t_{1}\}$ but no other terminal  is $d_{1}-m-(n-2m)$. Combining  $n\le d_{1}$ with all these numbers, we get that 
\[d_{1}-m-(n-2m)=d_{1}-n+m\ge d_{1}-d_{1}+m=m.\]
The claim is proved.
\end{claimproof}  

\begin{claim} \label{cl:fromB1toA1} Suppose that there are at most $d_1+2$ terminals in $B_1=K^{d_1}\times K^2$. If every row  in the subgraph $A_{1}=K^{d_{1}}\times K^{d_{2}-1}$ of $\bar C_{12}$ has a free entry,  then, for every terminal $x\not\in \{s_1,t_{1}\}$ in $B_{1}$,  there is a $B_{1}-A_{1}$ path $L$ from $x$ to a free entry in $A_{1}$ such that $L$ is $X$-valid; and all these $X$-valid paths are disjoint.    
\end{claim}
   

\begin{claimproof}
If a row of $B_{1}$ contains exactly one terminal $x\not\in \{s_1,t_{1}\}$, then send $x$ to a free entry in the same row of $A_{1}$. Let $x_{1}$ and $x_{2}$ be two terminals in $B_{1}$ that satisfy $\left\{x_{1},x_{2}\right\}\cap \left\{s_1,t_{1}\right\}=\emptyset$ and occupy a row $r_{f}$ of $B_{1}$. From \cref{cl:full-free-rows} ensues the existence of a row $r_{e}$ of $B_{1}$ that contains either no terminal or a terminal in $\{s_1,t_{1}\}$ but no other terminal; in short, there is at least a free entry in $r_{e}$. 

Consider a pair $(r_{f},r_{e})$ of rows granted by \cref{cl:full-free-rows}. Send either $x_{1}$ or $x_{2}$, say $x_{1}$, to the free entry in the row $r_{e}$ of $A_{1}$ passing through the corresponding free entry in the row $r_{e}$ of $B_{1}$, and send $x_{2}$ to a free entry in the row $r_{f}$ of $A_{1}$.  The proof of the claim is now complete.
\end{claimproof}

Now suppose that $\alpha=0$ or $2\le \alpha\le d_1-2$. In this subcase, the subgraph $\bar C_{12}$ contains at most $\alpha$ full rows: if $\alpha+1$ rows were full in $\bar C_{12}$ then there would be at least $(\alpha+1)(d_{2}-1)$ terminals in $\bar C_{12}$ but $(\alpha+1)(d_{2}-1)>d_{2}-2+\alpha$  (\cref{lem:inequality}). Even when the path $L_{1}$ uses the first row of $C_{12}$ by (\ref{eq:path-L1}), there is no loss of generality by assuming that the full rows of $\bar C_{12}$ are among  the {\it first $\alpha+1$ rows} of $\bar C_{12}$. It follows that every row of $A_{\alpha+1}$ has a free entry. 

We  send to $B_{\alpha+1}$ the terminals other than $s_1$ and $t_{1}$ that are in the rows 2 to $\alpha+1$ of $C_{12}$, the terminals other than $s_1$ and $t_{1}$ that are in  $B_{1}\setminus B_{\alpha+1}$.  For $\alpha=0$, $B_{1}\setminus B_{\alpha+1}=\emptyset$ and there is nothing to do.   

We now focus on the subcase $2\le \alpha\le d_{1}-2$. Let $n_{1}$ and $n_{2}$ denote the number of terminals in $B_{1}\setminus B_{\alpha+1}$ and $B_{\alpha+1}$, respectively. Then the following inequalities hold 
\begin{align*}
n_{1}+n_{2}&\le d_{1}+2-\alpha \le d_{1}\quad \text{(since $ 2\le \alpha$)},\\
n_{1}+n_{2}&\le d_{1}+2-\alpha\le 2d_{1}-2\alpha=|B_{\alpha+1}| \quad \text{(since $\alpha\le d_{1}-2$)}.
\end{align*}
Since  $B_{1}$ is $d_{1}$-connected by \cref{lem:connectivity-complete-graphs}, Menger's theorem gives $n_{1}$ disjoint paths from the terminals in $B_{1}\setminus B_{\alpha+1}$ to $n_{1}$ free entries in $B_{\alpha+1}$. These free entries in $B_{\alpha+1}$ will henceforth play the role of the terminals of $B_{1}\setminus B_{\alpha+1}$, and so we call them terminals as well. 

Let $d_{1}':=d_{1}-\alpha$ and $d_{2}':=d_{2}$. Then, $d_1'\ge 2$, there are at most $d_{1}'+2$ terminals in $B_{\alpha+1}=K^{d_{1}'}\times K^{2}$, and every row in $A_{\alpha+1}=K^{d_{1}'}\times K^{d_2'}-1$ has a free entry.  Hence, \cref{cl:fromB1toA1} applies, and there is a linkage formed by $X$-valid paths from the terminals in $B_{\alpha+1}$, other than $t_{1}$, to free entries in $A_{\alpha+1}$.  Now we  have a new linkage problem $Y'$ in $\bar C_{12}$ with at most $2(k-1)$ pairs. The  solution of $Y'$ in $\bar C_{12}$ implies a solution of the linkage problem $Y$ in $G$. To link the pairs of $Y'$  use the $(k-1)$-linkedness of $\bar C_{12}$.

Finally assume that $\alpha=1$. In a first scenario suppose  that either both entries in $B_1\setminus B_2$ are nonterminals or each  terminal in $B_1\setminus B_2$ is adjacent to a nonterminal in $B_{2}$. Then we can send the terminals in $B_1\setminus B_{2}$ to $B_{2}$. In the second scenario suppose that there is a terminal $s_{i}$ in $B_1\setminus B_2$ whose neighbours in  $B_{2}$ are all terminals. Then the column of $s_1$ in $B_{1}$ would contain exactly  $d_{1}$ terminals, including $s_{i}$. We send $s_i$ to a free entry in $A_1$, in the same row as $s_i$ (the first row of $A_1$): if this free entry didn't exist, then $s_{i}$ would be  adjacent to the $d_{2}-1$ terminals in $A_{1}$ and the $d_{1}-1$ terminals in $B_{2}$. Since there are $d_1+d_2$ terminals in total, it would follow that  $s_{i}$ is adjacent to $t_{i}$. This contradiction shows that we can send $s_i$ to a free entry in $A_1$. This completes the proof of \cref{case:d1+2}. 

In both scenarios, it remains to send the terminals in $B_2$ to $A_2$. To do so, we reason as in the subcase $2\le \alpha\le d_{1}-2$.  Let $d_{1}':=d_{1}-1$ and $d_{2}':=d_{2}$. It follows that $d_{1}'\ge 2$, that there are at most $d_{1}+2'$ terminals in $B_{2}$, and that every row in $A_{2}$ has a free entry.  \cref{cl:fromB1toA1} applies again and gives a linkage formed by $X$-valid paths from the terminals in $B_{2}$, other than $s_1,t_{1}$, to free entries in $A_{2}$.  

With all the terminals in $\bar C_{12}$, therein we have a new linkage problem $Y'$ with at most $2(k-1)$ pairs whose solution in $\bar C_{12}$ implies a solution of the linkage problem $Y$ in $G$. To solve  $Y'$ in $\bar C_{12}$ use the $(k-1)$-linkedness of $\bar C_{12}$. 

By symmetry, we also have the result if there are at most $d_{2}+2$  terminals in $R_{12}$, including $\{s_{1},t_{1}\}$.

\begin{case} The subgraph $C_{12}$ contains at least $d_{1}+3$ terminals, including $\{s_{1},t_{1}\}$. 
\end{case}

This case reduces to the previous case. If  $C_{12}$ contains at least $d_{1}+3$ terminals then $R_{12}$ contains at most $d_{2}-3+4=d_{2}+1$ terminals, since there are four entries shared by $C_{12}$ and $R_{12}$. Since we make no distinction between columns and rows, this case is already covered. This completes the proof of the theorem.
\end{proof}

\section{Duals of cyclic polytopes}

There is a close connection between duals of cyclic $d$-polytopes with $d+2$ vertices  and Cartesian products of complete graphs.

The {\it moment curve} in $\mathbb{R}^d$ is defined by $x(t):=(t,t^{2},\ldots,t^{d})$ for $t\in \mathbb{R}$, and the convex hull of any $n>d$ points on it gives a {\it cyclic polytope} $C(n,d)$. The {\it combinatorics} of a cyclic polytope, the face lattice of the polytope faces partially ordered by inclusion, is independent of the points chosen on the moment curve. Hence we talk of the cyclic $d$-polytope on $n$ vertices \citep[Example 0.6]{Zie95}.

For a polytope $P$ that contains the origin in its interior, the \textit{dual polytope} $P^*$ is  defined as 
\[
P^*=\{y\in \R^d\mid x\cdot y\le 1 \text{ for all $x$ in $P$}\}.
\] 
If $P$ does not contain the origin, we translate the polytope so that it does. Translating the polytope $P$ changes the geometry of $P^*$ but not its face lattice.  The face lattice of $P^*$ is the inclusion reversed face lattice of $P$. In particular, the vertices of $P^*$ correspond to the facets of $P$, and the edges of $P^*$ correspond to the $(d-2)$-faces of $P$. The \emph{dual graph} of a polytope $P$ is the graph of the dual polytope, or equivalently, the graph on the set of facets of $P$ where two facets are adjacent in the dual graph if they share a $(d-2)$-face. 

Duals of cyclic $d$-polytopes are simple $d$-polytopes. It is also the case that the dual of a  cyclic $d$-polytope with $d+2$ vertices can be expressed as $T({\floor{d/2}})\times T({\ceil{d/2}})$ ({\citealt[Ex.~0.6]{Zie95}}).  From this observation and \cref{thm:simplex-products} the next corollary follows at once.

\begin{corollary}
\label{cor:linkedness-cyclic-d+2}  Duals of cyclic polytopes with $d+2$ vertices are $\floor{d/2}$-linked for every $d\ge 2$.
\end{corollary}

\section{Acknowledgments}

Guillermo would like to thank the hospitality of Leif J{\o}rgensen, and the Department of Mathematical Sciences at Aalborg University where this research started. 

Julien Ugon's research was supported by the ARC discovery project DP180100602. 
 
\bibliographystyle{apalike2}

\end{document}